\documentclass[a4paper]{amsart}
\usepackage{amssymb,amscd}
\usepackage{fullpage}

\newtheorem{thm}{Theorem}[section]
\newtheorem{prop}[thm]{Proposition}
\newtheorem{lem}[thm]{{Lemma}}
\newtheorem{defn}[thm]{Definition}

\newtheorem{cor}[thm]{Corollary}
\newtheorem{conj}{Conjecture}
\numberwithin{equation}{section}

\def\.{\cdot}
\def\<{\left\langle}
\def\>{\right\rangle}
\def\({\left(}
\def\){\right)}
\def\lra{\longrightarrow}
\renewcommand{\leq}{\leqslant}
\renewcommand{\geq}{\geqslant}
\renewcommand{\phi}{\varphi}

\def\bar#1{\overline{#1}}

\def\L{\mathbb\L}

\def\subset{\subseteq}

\def\epsilon{\varepsilon}
\def\U{\mathcal U}
\def\odd{\mathrm{odd}}

\def\ker{\operatorname{ker}}

\def\index{\operatorname{index}}

\def\im{\operatorname{im}}

\def\C{\mathbb C}
\def\D{\mathbb D}
\def\Z{\mathbb Z}

\def\R{\mathbb R}

\def\T{\mathbb T}
\def\S{\mathbb S}
\def\ds{\displaystyle}

\def\virt-dim{\operatorname{virt-dim}}

\def\Gr{\operatorname{Gr}}

\def\Id{\operatorname{Id}}

\def\ev{\operatorname{ev}}

\def\dim{\operatorname{dim}}

\title[Complex Cobordism of Hilbert Manifolds]
{On Fredholm Index,Transversal Approximations and Quillen's
Geometric Complex Cobordism of Hilbert Manifolds with some
Applications to Flag Varieties of Loop Groups}
\author{Cenap \"{O}zel}
\address{AIBU Golkoy Kampusu, Bolu 14280, Turkey.}
\email{cenap@ibu.edu.tr}
\thanks{We are indebted to Unit of Research and Developing,
Bolu Abant Izzet Baysal University for their financial support
during our project with numbered 2000/03/04/67. I would like to
thank A.J. Baker and Dept. of Math, University of Glasgow for personal support,conservations and hospitality.}%endthanks
%endkeywords
\subjclass{Algebraic topology,Lie theory,Global Analysis}%endsubjclass
\date{April 16,2002}

\begin{document}
\begin{abstract}
In \cite{baker-ozel}, by using Fredholm index we developed a
version of Quillen's geometric cobordism theory for infinite
dimensional Hilbert manifolds. This cobordism theory has a graded
group structure under topological union operation and has
push-forward maps for complex orientable Fredholm maps. In this
work, by using Quinn's Transversality Theorem \cite{Quinn}, it
will be shown that this cobordism theory has a graded ring
structure under transversal intersection operation and has
pull-back maps for smooth maps. It will be shown that the Thom
isomorphism in this theory will be satisfied for finite
dimensional vector bundles over separable Hilbert manifolds and
the projection formula for Gysin maps will be proved. After we
discuss the relation between this theory and classical cobordism,
we describe some applications to the complex cobordism of flag
varieties of loop groups and we do some calculations.
\end{abstract}
\maketitle
\section{Preliminaries.}
%Fredholm Maps and Cobordism of  Separable Hilbert Manifolds

\subsection{Complex Cobordism.}
Complex bordism theory $MU$ was originally defined by geometric
means as bordism classes of maps of stably complex manifolds. For
a space $X$, $MU_q(X)$ is the set of equivalence classes of maps
$M\stackrel{f}\rightarrow X$ where $M$ is a closed, stably almost
complex manifold. It means that $M$ is a compact smooth manifold
without boundary of dimension $q$ with $TM$ stably complex.

Two such maps $(M,f)$ and $(N,g)$ are bordant iff their
topological union extends to a map $W\rightarrow X$ of a compact
stably almost complex manifold $W$ of dimension $q + 1$, whose
boundary is the union of $M$ and $N$. Here, the stably almost
structures on $M$ and $N$ induced by the embedding into $W$ and
the original ones are required to be equivalent. The Abelian group
structure on $MU_* (X)$ is given by the operation of disjoint
union. In \cite{thom}, Rene Thom gave the homotopy theoretic
construction of this group. The details can be found in
\cite{stong}.

The dual cohomology theory $MU$ called complex cobordism was a
given a geometric description by D. Quillen. His detailed
construction can be found in \cite{Quillen}. For a manifold $X$ of
dimension $n$ an element in $MU^{n-q}(X)$ is represented by a
smooth proper map $M\rightarrow X$ of a (not necessiraly compact)
manifold $M$ of dimension $q$ together with an equivalence class
of complex orientations. Two such maps are cobordant iff they are
bordant as maps with complex orientations.

Now we recall the basic properties of multiplicative generalized
cohomology theories with complex orientation.

A multiplicative cohomology theory $h$ is a functor from
topological spaces to graded rings satisfying the
Eilenberg-Steenrod axioms. Details can be found in \cite{Dyer} and
\cite{adams}. A multiplicative cohomology theory $h$ is complex
oriented if the complex vector bundles are oriented for $h$. In a
complex oriented cohomology theory the Euler class and the Chern
classes of complex vector bundles are defined and satisfy the
usual properties. Examples of complex oriented cohomolgy theories
are ordinary cohomology, complex $K$-theory, elliptic cohomology
and complex cobordism.

Complex cobordism is the universal complex oriented cohomology
theory and so for any such theory $h$ we have canonical map
$MU\rightarrow h$.

A complex orientation of a proper map of smooth manifolds $f:
M\rightarrow N$ is a factorization
$$
M  \stackrel{i}\rightarrow  \xi \stackrel{\pi}\rightarrow  N
$$
where $\xi$ is a complex vector bundle and $i$ is an embedding
with a stably complex normal bundle.

A compact manifold is said to be complex oriented iff the tangent
bundle is stably complex.

All details about complex cobordism and multiplicative complex
oriented generalized cohomology theories can be found in
\cite{adams},\cite{Dold},\cite{Dyer} and \cite{Quillen}.

\subsection{Infinite Dimensional Manifolds and Pressley-Segal Stratifications}
The general reference for the section is \cite{segal-pressley}.

$LG$ is the group of all smooth maps from $\S^1$ to compact
semi-simple Lie group $G$. Its essential homogeneous spaces are
infinite flag manifold $LG/T$ and based loop space $\Omega G$
where $T$ is a maximal torus of $G$. We are interested to infinite
flag manifold $LG/T$. The cell complexes and the dual
stratifications of $LG/T$ are analogically like as $G/T$. Also we
are interested to Grassmannians $\Gr (H)$ of infinite dimensional
complex separable Hilbert space $H$. Its stratifications can be
found in \cite{segal-pressley}.
\section{The Fredholm Index and Complex Cobordism of Hilbert Manifolds.}
In \cite{Quillen}, Quillen gave a geometric interpretation of
cobordism groups which suggests a way of defining the cobordism of
separable Hilbert manifolds equipped with suitable structure. In
order that such a definition be sensible, it ought to reduce to
his for finite dimensional manifolds and smooth maps of manifolds
and be capable of supporting reasonable calculations for important
types of infinite dimensional manifolds such as homogeneous spaces
of free loop groups of finite dimensional Lie groups.

\subsection{Cobordism of separable Hilbert  manifolds.} By a
manifold, we mean a smooth manifold modelled on a separable
Hilbert space; see Lang \cite{lang} for details on infinite
dimensional manifolds. The facts about Fredholm map can be found
in \cite{conway}.
\begin{defn}\label{zirto}
Suppose that $f: X\rightarrow Y$ is a proper Fredholm map with
even index at each point. Then $f$ is an \emph{admissible complex
orientable map} if there is a smooth factorization
$$f: X\stackrel{\tilde{f}}\rightarrow \xi \stackrel{q}\rightarrow Y,$$
where $q:\xi \rightarrow Y$ is a finite dimensional smooth complex
vector bundle and $\tilde{f}$ is a smooth embedding endowed with a
complex structure on its normal bundle $\nu (\tilde{f})$.
\par
A complex orientation for a Fredholm  map $f$ of odd index is
defined to be one for the map $ (f,\epsilon): X\rightarrow Y\times
\mathbb R $ given by  $(f,\epsilon) (x) = (f(x),0)$ for every $x
\in X$. At $x\in  X$, $\index (f, \epsilon)_x = (\index f_x) - 1$.
Also the finite dimensional  complex vector bundle $\xi$ in the
smooth factorization will be replaced by $\xi \times \mathbb R$.
\end{defn}
Suppose that $f$ is an admissible complex orientable map. Then
since the map $f$ is the Fredholm and $\xi$ is a finite
dimensional vector bundle, we see $\tilde{f}$ is also a Fredholm
map. By the surjectivity of $q$,
$$\index \tilde{f} = \index f - \dim \xi.$$
Before we give a notion of equivalence of such factorizations
$\tilde{f}$ of $f$, we want to give some definitions.
\begin{defn}
Let $X$, $Y$ be the smooth separable Hilbert manifolds and $F: X
\times \mathbb R \rightarrow Y$ a smooth map. Then we will say
that $F$ is an \emph{isotopy} if it satisfies the following
conditions.
\begin{enumerate}
\item For every $t \in \mathbb R$, the map $F_t$ given by
$F_t (x) = F(x, t)$ is an embedding.
\item There exist numbers
$t_0 < t_1$ such that $F_t = F_{t_0}$ for all $t \leq t_0$ and
$F_t = F_{t_1}$ for all $t \geq t_1$.
\end{enumerate}
The closed interval $[t_0, t_1]$ is called a \emph{proper domain}
for the isotopy. We say that two embeddings $f: X \rightarrow Y$
and $g: X \rightarrow Y$ are \emph{isotopic} if there exists an
isotopy $F_t: X\times \mathbb R \rightarrow Y$ with proper domain
$[t_0, t_1]$ such that $f = F_{t_0}$ and $g = F_{t_1}$.
\end{defn}
\begin{prop}(see \cite{lang})\label{lang}
The relation of isotopy between smooth embeddings is an
equivalence relation.
\end{prop}
\begin{defn}
Two factorizations $f: X\stackrel{\tilde{f}}\rightarrow \xi
\stackrel{q}\rightarrow Y$ and $f:
X\stackrel{\tilde{f'}}\rightarrow \xi' \stackrel{q'}\rightarrow Y$
are \emph{equivalent} if $\xi$ and $\xi'$ can be embedded as
subvector bundles of a vector bundle $\xi''\rightarrow Y$ such
that $\tilde{f}$ and $\tilde{f'}$ are isotopic in $\xi''$ and this
isotopy is compatible with the complex structure on the normal
bundle. That is, there is an isotopy $F$ such that for all $t \in
[t_0, t_1]$, $F_t:  X \rightarrow \xi''$ is endowed with a complex
structure on its normal bundle which matches that of $\tilde{f}$
and $\tilde{f'}$ in $\xi''$ at $t_0$ and $t_1$ respectively.
\end{defn}
By Proposition \ref{lang}, we have
\begin{prop}
The relation of equivalence of admissible complex orientability of
proper Fredholm maps between separable Hilbert manifolds is an
equivalence relation.
\end{prop}
This generalizes  Quillen's notion of complex orientability for
maps of finite dimensional manifolds. We can also define a notion
of cobordism of admissible complex orientable maps between
separable Hilbert manifolds. First we recall some ideas on the
transversality.

\begin{defn}Let $f_1: M_1 \rightarrow N, f_2: M_2 \rightarrow N$
be smooth maps between Hilbert manifolds. Then $f_1$ and $f_2$ are
\emph{transverse} at $y \in N$ if
\[
df_1 (T_{x_1} M_1) + df_2 (T_{x_2}M_2) = T_y N
\]
whenever $f_1 (x_1) = f_2(x_2) = y$. The maps $f_1$ and $f_2$ are
said to be \emph{transverse} if they are transverse at every point
of $N$.
\end{defn}
\begin{lem}
Smooth maps $f_i: M_i \rightarrow N (i = 1,2)$ are transverse if
and only if $f_1\times f_2 : M_1 \times M_2 \rightarrow N\times N$
is transverse to the diagonal map $\Delta :  N\rightarrow N\times
N$.
\end{lem}
\begin{defn}
Let $f_1: M_1 \rightarrow N, f_2: M_2 \rightarrow N$ be transverse
smooth maps between smooth Hilbert manifolds. The
\emph{topological pullback}
\[
M_1 \prod_N M_2 = \{ (x_1 , x_2) \in M_1 \times M_2 : f_1(x_1) =
f_2(x_2)\}
\]
is a submanifold of $M_1 \times M_2$ and the diagram
\[
\begin{CD}
M_1 \prod_N M_2  @>{f_2}^*(f_1)>> M_2  \\
@VV{f_1}^*(f_2)V         @VVf_2V \\
M_1 @>f_1>> N
\end{CD}
\]
is commutative, where the map ${f_i}^* (f_j)$ is pull-back of
$f_j$ by $f_i$.
\end{defn}

\begin{defn}
Let $f_i: X_i \rightarrow  Y (i= 0,1)$ be admissible complex
oriented maps. Then $f_0$ is \emph{cobordant} to $f_1$ if there is
an admissible complex orientable map $h: W \rightarrow Y\times
\mathbb R$ such that the maps $\epsilon_i: Y\rightarrow Y \times
\mathbb R \,\text{given by} \,\epsilon_i (y) = (y,i)$ for $i= 0, 1
$, are transverse to $h$ and the pull-back map ${\epsilon_i}^* h$
is equivalent to $f_i$. The cobordism class of $f: X\rightarrow Y$
will be denoted by $[X, f]$.
\end{defn}
\begin{prop}\label{grem}
If $f: X\rightarrow Y$ is an admissible complex orientable map and
$g: Z\rightarrow Y$ a smooth map transverse to $f$, then the
pull-back map
\[
{g^*}(f):  Z\prod_Y X\rightarrow Z
\]
is an admissible complex orientable map with finite dimensional
pull-back vector bundle
\[
g^*(\xi) = Z\prod_Y \xi = \{(z,v)\in Z\times \xi : g(z) = q(v)\}
\]
in the factorization of $g^* (f)$, where $q: \xi\rightarrow Y$ is
the finite-dimensional complex vector bundle in the factorization
of $f$ as in Definition \ref{zirto}.
\end{prop}
The next result was proved in \cite{ozel} by essentially the same
argument as in the finite dimensional situation using the Implicit
Function Theorem \cite{lang}.
\begin{thm}
Cobordism is an equivalence relation.
\end{thm}

\begin{defn}
For a separable Hilbert manifold $Y$, $\mathcal U^d (Y)$ is the
set of cobordism classes of  the admissible complex orientable
proper Fredholm maps of index $-d$.
\end{defn}

In the above definition, instead of proper maps, closed maps could
be used for infinite dimensional Hilbert manifolds, because of the
following result of Smale \cite{smale}.
\begin{thm}\label{martin}
When $X$ and $Y$ are infinite dimensional, every closed Fredholm
map $X\rightarrow Y$ is proper.
\end{thm}
My next result is the following.

\begin{thm}\label{iain}
If $f: X\rightarrow Y$ is an admissible complex orientable
Fredholm map of index $d_1$ and $g: Y\rightarrow Z$ is an
admissible complex orientable Fredholm map of index $d_2$, then
$g\circ f: X\rightarrow Z$ is an admissible complex orientable
Fredholm map with index $d_1 + d_2$.
\end{thm}

Let $g:Y \rightarrow Z$ be an admissible complex orientable
Fredholm map of index $r$. By Theorem \ref{iain}, we have
\emph{push-forward, or Gysin map}
\[
g_{*}: \mathcal U^d (Y) \rightarrow \mathcal U^{d + r} (Z)
\]
given by $g_{*} ([X, f]) = ([X, g \circ f])$.

We show in \cite{ozel} that it is well-defined. If $g':
Y\rightarrow Z$ is a second map cobordant to $g$ then ${g'}_* =
g_*$; in particular, if $g$ and $g'$ are homotopic through proper
Fredholm maps they induce the same Gysin maps. Clearly, we have
$(h\circ g)_{*} = h_{*}g_{*}$ for admissible complex orientable
Fredholm maps $h,g$ and $\Id_{*} = \Id$.

The graded cobordism set $\mathcal U^* (Y)$ of the separable
Hilbert manifold $Y$ has a group structure given as follows. Let
$[X_1 , f_1]$ and $[X_2 , f_2]$ be cobordism classes. Then $[X_1 ,
f_1] + [X_2 , f_2]$ is the class of the map $f_1 \sqcup f_2: X_1
\sqcup X_2 \rightarrow Y$, where $X_1 \sqcup X_2$ is the
topological sum (disjoint union) of $X_1$ and $X_2$.  We show in
\cite{ozel} that this sum is well-defined. As usual, the class of
the empty set $\emptyset$ is the zero  element of  the cobordism
set and the negative of $[X,f]$ is itself with the opposite
orientation on the normal bundle of the embedding $\tilde{f}$.
Then we have
\begin{thm}
The graded cobordism set $\mathcal U^* (Y)$ of the admissible
complex orientable maps of $Y$ is a graded abelian group.
\end{thm}
Now we define relative cobordism .
\begin{defn}
If $A$ is a finite dimensional submanifold of $Y$, the relative
cobordism set $\mathcal U^* (Y,A)$ is the set of the admissible
complex orientable maps of $Y$ whose images lie in $Y-A$.
\end{defn}
More generally,
\begin{thm}\label{relative}
Let $A$ be a finite dimensional submanifold of $Y$. Then the
relative cobordism set $\mathcal U^* (Y,A)$ is a graded abelian
group and there is a homomorphism $\kappa^* :\mathcal U^*
(Y,A)\rightarrow \mathcal U^* (Y)$ by
$\kappa^*[M\stackrel{h}\rightarrow Y] =[M\stackrel{h}\rightarrow
Y]$ with $h(M)\subset Y-A$.
\end{thm}

If our cobordism  functor $\mathcal U^*( \,)$ of admissible
complex orientable Fredholm maps is restricted to finite
dimensional Hilbert manifolds, it agrees Quillen's complex
cobordism functor $MU^* (\,)$.
\begin{thm}
For finite dimensional separable Hilbert manifolds $A\subset Y$,
there is a natural isomorphism
$$
\mathcal U^*(Y,A)\cong MU^*(Y,A).
$$
\end{thm}
\section{Transversal approximations, contravariance and cup products.} We
would like to define a product structure on the graded cobordism
group $\mathcal U^* (Y)$. Given  cobordism classes $[X_1 , f_1]
\in \mathcal U^{d_1}(Y_1)$ and $[X_2 , f_2]\in \mathcal
U^{d_2}(Y_2)$, their external product is
\[
[X_1 , f_1]\times [X_2 , f_2 ] = [X_1 \times X_2 ,f_1 \times f_2]
\in \mathcal U^{d_1 + d_2} (Y_1\times Y_2).
\]
Although there is the external product in the category of
cobordism of separable Hilbert manifolds, we can not necessarily
define an internal product on $U^{*}(Y)$ unless $Y$ is a finite
dimensional manifold. However,if admissible complex orientable
Fredholm map $f_1 \times f_2$ is transverse to the diagonal
imbedding $\Delta:Y\rightarrow Y\times Y$, then we do have an
internal (cup) product
\[
[X_1 , f_1] \cup [X_2 , f_2] = \Delta^* [X_1\times X_2, f_1\times
f_2].
\]
If $Y$ is finite dimensional, then by Thom's Transversality
Theorem in \cite{thom}, every complex orientable map to $Y$ has a
transverse approximation, hence the cup product $\cup$ induces a
graded ring structure on $\mathcal U^* (Y)$. The unit element $1$
is represented by the identity map $Y\rightarrow Y$ with index
$0$. However F. Quinn \cite{Quinn} proved the generalization of
Thom's Transversality Theorem for separable Hilbert manifolds
using smooth transversal approximations of Sard functions in fine
topology.

Details about the fine topology, jets  and smooth maps space
$C^m(W, N)$,$C^{\infty}(W, N)$ can be found in \cite{michor}. In
this topology, the derivatives of the difference function between
the function $g$ and its approximation $g'$ are bounded. We would
like to interpret this approximation in the fine topology. We need
some notation to describe this situation.

\begin{defn}
Let $X$ and $Y$ be smooth manifolds. A \emph{k-jet from} $X$
\emph{to} $Y$ is an equivalence class $[f, x]_k$ of pairs $(f, x)$
where $f: X\rightarrow Y$ is a smooth mapping, $x \in X$. The
pairs $(f, x)$ and $(f', x')$ are \emph{equivalent} if $x = x'$,
$f$ and $f'$ have same Taylor expansion of order $k$ at $x$ in
some pair of coordinate charts centered at $x$ and $f (x)$
respectively. We will write $J^k f (x)=[f, x]_k $ and call this
the \emph{k-jet of} $f$ at $x$.
\end{defn}
There is an equivalent definition of this equivalence relation:
$[f, x]_k = [f', x']_k$ if $x=x'$ and $T^k_x f= T^k_x f'$ where
$T^k$ is the $k$th tangent mapping.
\begin{defn}
For a topological space $X$,a covering of $X$ is \emph{locally
finite} if every point has a neighborhood which intersects only
finitely many elements of the covering.
\end{defn}
Approximation $g'$ of $g$ in the smooth fine topology means the
following. Let $\{L_i\}_{i \in I}$ be a locally finite cover of
$W$. For every open set $L_i$, there is a bounded continuous map
$\varepsilon_i: L_i \rightarrow [0, \infty)$ such that for every
$x \in L_i$ and $k > 0$,
\[
||J^k g (x) - J^k g' (x) || < \varepsilon_i (x).
\]

\begin{defn}
Let $E$ be a Banach space. We say that a collection $\mathcal S$
of smooth functions $\alpha: E \rightarrow \mathbb R$ is a
\emph{Sard class} if it satisfies the following conditions:
\begin{enumerate}
\item for $r \in \mathbb R$, $y \in E$ and $\alpha \in \mathcal S$,
then the function $x \rightarrow \alpha(rx + y)$ is also in the
class $\mathcal S$,
\item if $\alpha_n \in \mathcal S$, then the rank of differential
$D_x (\alpha_1, \ldots, \alpha_n)$ is constant for all $x$ not in
some closed finite dimensional submanifold of $E$.
\end{enumerate}
\end{defn}
\begin{defn}\label{Sard}
Let $\mathcal S$ be a Sard class on $E$, $U$open in $E$, and $M$ a
smooth Banach manifold. We define $\mathcal S(U, M)$ to be the
collection of \emph{Sard functions} $f: U\rightarrow M$ such that
for each $x \in U$ there is a neighbourhood  $V$ of $x$, functions
$\alpha_1, \ldots, \alpha_n \in \mathcal S$, and a smooth map
$g:W\rightarrow M$, where $W$ open in $\mathbb R^n$ contains
$(\alpha_1, \ldots, \alpha_n)(V)$, all such that $f|V= g \circ
(\alpha_1, \ldots, \alpha_n)|V$.
\end{defn}
\begin{defn}
The \emph{support} of a function $f: X\rightarrow \mathbb R$ is
the closure of the set of points $x$ such that $f(x)\neq 0$.
\end{defn}
From \cite{Quinn}, we have
\begin{thm}
$E$ admits a Sard class $\mathcal S$ if $\mathcal S(E, \mathbb R)$
contains a function with bounded nonempty support. In particular,
the separable Hilbert space admits Sard classes.
\end{thm}
\begin{defn}
A \emph{refinement} of a covering of $X$ is a second covering,
each element of which is contained in an element of the first
covering.
\end{defn}
\begin{defn}
A topological space is \emph{paracompact} if it is Hausdorff, and
every open covering has a locally finite open refirement.
\end{defn}
\begin{defn}
A smooth \emph{partition of unity} on a manifold $X$ consists of
an covering $\{U_i\}$ of $X$ and a system of smooth functions
$\psi_i: X\rightarrow \mathbb R$ satisfying the following
conditions.
\begin{enumerate}
\item $\forall x \in X$, we have $\psi_i (x) \geq 0$;\\
\item  the support of $\psi_i$ is contained in $U_i$;\\
\item the covering is locally finite;\\
\item for each point $x \in X$, we have
\[
\ds\sum_{i} \psi_i (x) = 1.
\]
\end{enumerate}
\end{defn}
\begin{defn}
A manifold $X$ will be said to \emph{admit partitions of unity} if
it is paracompact, and if, given a locally finite open covering
$\{U_i\}$, there exists a partition of unity $\{\psi_i\}$ such
that the support of $\psi_i$ is contained in some $U_i$.
\end{defn}
From \cite{lang}, we have
\begin{thm}
Let $X$ be a paracompact smooth manifold modelled on a separable
Hilbert space $H$. Then $X$ admits smooth partitions of unity.
\end{thm}
From \cite{Eells-McAlpin},
\begin{thm}
On a separable Hilbert manifold the functions constructed using
the partitions of unity form a Sard class.
\end{thm}

The following result was proved by F. Quinn .
\begin{thm}\label{kirmizi}
Let $H$ be the smooth separable Hilbert space and let $U$ be an
open set in $H$. If $f: W\rightarrow N$ is a smooth proper
Fredholm map, then smooth maps transversal to $f$ are dense in
$\mathcal S(U, N)$ with the $C^0$ fine topology.
\end{thm}
We will require the Open Embedding Theorem of Eells $\&$ Elworthy
\cite{Eells-Elw}.
\begin{thm}\label{open embedding}
Let $M$ be a smooth manifold modelled on the separable infinite
dimensional Hilbert space $H$. Then $M$ is diffeomorphic to an
open subset of $H$.
\end{thm}
\begin{thm}
The space $\mathcal S(U, N)$ is dense in $C^0 (U, N)$ in the $C^0$
fine topology.
\end{thm}
Using these techniques, Quinn proves the following result.
\begin{cor}\label{canalici}
Let $M$ be a smooth separable Hilbert manifold. If $f:
W\rightarrow N$ is a smooth proper Fredholm map, then smooth maps
transversal to $f$ are dense in $C^0(M,N)$ in the $C^0$ fine
topology.
\end{cor}

From \cite{Eells-Elw}, we have
\begin{thm}\label{homotopy}
Let $X$ and $Y$ be two smooth manifold modelled on the separable
infinite dimensional Hilbert space $H$. If there is a homotopy
equivalence $\phi: X\rightarrow Y$, then $\phi$ is homotopic to a
diffeomorphism.
\end{thm}
Now I will try to tell a detailed explanation of Quinn's very
technical work in Theorem \ref{kirmizi}.

Let $f: W\rightarrow N$ be a smooth proper Fredholm map and let $g
\in \mathcal S(U, N)$ specified in Theorem \ref{kirmizi}. Now we
are given an arbitrary $C^0$ fine neighborhood of $g$ in which we
want to find a map transversal to $f$. Since some neighborhood of
the image of $g$ is metrizable, we can impose a metric on it and
then we can choose a function $\varepsilon: U \rightarrow (0,1)$
so that the $\varepsilon$-neighborhood of $g$ in $\mathcal S(U,
N)$ lies in the given $C^0$ fine neighborhood. Now we explain the
Smale decomposition of the smooth proper Fredholm map $f:
W\rightarrow N$. If $x \in N$, then there is a coordinate
neighborhood $\Theta : U\thickapprox H \times \R^k$ about $x$ such
that
$$
f^{-1}(u)\supseteqq \ds\bigcup_{i=1}^n V_i \supseteqq
f^{-1}({\Theta}^{-1}(H \times \D^k))
$$
for some sets $V_i$ such that $\Psi_i : V_i \thickapprox H\times
W_i, W_i$ open in $\R^m$, and $\Theta \circ f \circ \Psi_{i}^{-1}
= (\pi,f_i): H\times W_i \rightarrow H\times \R^k$, where $\D^k$
is the open unit ball in $\R^k$.

Using separability of $U$, we cover $g(u)$ by $N-f(W)$ and a
countable number of sets of the form $\Theta_{i}^{-1}(H_i \times
\D^{k_i})$, where $\Theta_i : u_i \thickapprox H_i \times
\R^{k_i}$ are coordinate neighborhoods  as given in the Smale
decomposition of $f$. We denote the corresponding coordinate
neighborhoods in the domain by $(V_{i_j}, \Psi_{i_j})$,where
$\Psi_{i_j}:V_{i_j}\thickapprox H_i \times W_i$, $W_i$ open in
$\R^{m_i}$. Let $\{Y_i\}$ be a locally finite refinement of
$\{g^{-1}(\Theta_{i}^{-1}(H_i \times \D^{k_i})),
g^{-1}(N-f(W))\}$, and let $\{Z_i\}$ be a subcover such that $Y_i
\supseteqq \bar{Z_i}$.

Inductively we will get an approximation of $g$. Let $g_0 = g$.
Given $g_i$, let $g_{i+1}$ be the approximation defined by Quinn
applied the situation $U_1 =U_2 = Y_i, U_3= Z_i$, $W$ the disjoint
union of of the $W_{i_j}$ over $j$, and the $C^0$ fine topology
neighborhood so small that $g_{i+1} = g_i +
\frac{\varepsilon}{2^{i+1}}y_i$, and $g_{i+1}(Y_i) \subset
\Theta_{i}^{-1}(H_i \times \D^{k_i})$, where $y_i \in \S^{k_i}$.
Since $\{Y_i\}$ is locally finite, $g' =
\ds\lim_{i\rightarrow\infty}g_i$ is well-defined. $g'$ is an
$\varepsilon$-approximation of $g$ and $g'\in\mathcal S (U, N)$
and it is transverse to $f$ everywhere. It is interesting that
this approximation can be done in the $C^r$ fine topology. For
separable Hilbert manifolds, it can be even done in the smooth
$C^{\infty}$ fine topology. In this case, $g' \in
\overline{\mathcal S(U, N)}$.
\begin{thm}\label{hist}
Let $U$ be an open set in separable infinite dimensional Hilbert
space $H$ and let $f: M\rightarrow N$ be a proper Fredholm map
between separable infinite dimensional Hilbert manifolds $M$ and
$N$. Then the set of maps transverse to $f$ is dense in the
closure of Sard function space $\overline{\mathcal S(U, N)}$ in
the $C^{\infty}$ fine topology.
\end{thm}

By Corollary \ref{canalici}, a smooth map (even continuous map)
$g: Z\rightarrow Y$ can be deformed to a smooth map $g':
Z\rightarrow Y$ by a small correction until it is transverse to an
admissible complex orientable map $f: X \rightarrow Y$. It is
obvious that they are homotopic each other. By definition of
Cobordism and Proposition \ref{grem}, the cobordism functor is
contravariant for any smooth map between separable Hilbert
manifolds.
\begin{thm}\label{son}
Let $f: X\rightarrow Y$ be an admissible complex oriented map and
let $g: Z\rightarrow Y$ be a smooth (may be continuous) map. Then
the cobordism class of the pull-back $Z \prod_Y X \rightarrow Z$
depends only on the cobordism class of $f$, hence there is a map
$g^*: \mathcal U^d (Y)\rightarrow \mathcal U^d (Z)$ given by
$$g^* [X, f] =g'^* [X, f] = [Z\prod_Y X, {g'}^*(f)],$$ where $g'$ is
a smooth $\varepsilon$-approximation of g which is transverse to
$f$. Moreover, $g^*$ depends only on the homotopy class of $g$.
\end{thm}
Now we give the functorial property of $\mathcal U^*$ theory.
\begin{thm}
Let $X,Y$ and $Z$ be separable Hilbert manifolds. If
$Z\stackrel{\alpha} \rightarrow Y \stackrel{\beta}\rightarrow X$
are smooth functions, then
$$
(\beta \circ \alpha)^* = {\alpha}^* {\beta}^*: \mathcal
U^d(X)\rightarrow \mathcal U^d(Z).
$$
The identity map $\Id : X\rightarrow X$ induces the identity
endomorphism $\Id^* :\mathcal U^d(X)\rightarrow \mathcal U^d(X)$
for every $d$.
\end{thm}
\begin{proof}
After making them transverse where appropriate we consider the
following commutative diagram
\[
\begin{CD}
M''  @>{\alpha}'>> M' @>{\beta}'>> M \\
@VVf''V         @VVf'V @VVfV \\
Z @>\alpha >> Y @>\beta >> X,
\end{CD}
\]
where both squares are pullback diagrams with transverse $(\beta,
f)$ and $(\alpha, f')$; the outer square is then also pullback
with transverse $(\beta \circ\alpha, f)$. Given an admissible
complex oriented map $f: M\rightarrow X$, we show that the maps
$M'\stackrel{f'}\rightarrow Y$ and $M''\stackrel{f''}\rightarrow
Z$ are admissible complex orientable and  $\alpha^* \beta^*
[M\stackrel{f}\rightarrow X] =(\beta \alpha
)^*[M\stackrel{f}\rightarrow X]$. It is clear that $M''=
Z\ds\prod_Y M' = Z\ds\prod_Y (Y \ds\prod_X M)$ as well as
$Z\ds\prod_X M$ so that $\alpha^* \beta^*
[M\stackrel{f}\rightarrow X] =(\beta \alpha
)^*[M\stackrel{f}\rightarrow X]$ by the following commutative
diagram
\[
\begin{CD}
M''  @>{\alpha}'>> M' @>{\beta}'>> M \\
@VV\tilde{f''}V         @VV\tilde{f'}V @VV\tilde{f}V \\
\xi'' @>\alpha >> \xi' @>\beta >> \xi\\
@VVq''V         @VVq'V @VVqV \\
Z @>\alpha >> Y @>\beta >> X,
\end{CD}
\]
where the finite dimensional vector bundle $\xi''= Z\ds\prod_Y
\xi' = Z\ds\prod_Y (Y \ds\prod_X \xi)$ as well as $Z\ds\prod_X
\xi$.
\end{proof}

Let turn back the interior(cup) products in $\mathcal U^*$. Given
cobordism classes $[X_1,f_1] \in \mathcal U^{d_1}(Y_1)$ and $[X_2
, f_2]\in \mathcal U^{d_2}(Y_2)$, their external product is
\[
[X_1 , f_1]\times [X_2 , f_2 ] = [X_1 \times X_2 ,f_1 \times f_2]
\in \mathcal U^{d_1 + d_2} (Y_1\times Y_2).
\]
If admissible complex orientable Fredholm map $f_1 \times f_2$ is
transverse to the diagonal imbedding $\Delta:Y\rightarrow Y\times
Y$, then we do have an internal (cup) product
\[
[X_1 , f_1] \cup [X_2 , f_2] = \Delta^* [X_1\times X_2, f_1\times
f_2].
\]
If the diagonal imbedding $\Delta:Y\rightarrow Y\times Y$ is not
transverse to smooth proper Fredholm map $f_1 \times f_2:X_1
\times X_2 \rightarrow $, by Quinn's transversality Theorem, we
can find a smooth $\varepsilon$-approximation $\Delta'$ of
$\Delta$ which is transverse to $f_1\times f_2$. Then
\begin{thm}
If $[X_1,f_1] \in \mathcal U^{d_1}(Y_1)$ and $[X_2 , f_2]\in
\mathcal U^{d_2}(Y_2)$,internal(cup) product
\[
[X_1 , f_1] \cup [X_2 , f_2] = \Delta^* [X_1\times X_2, f_1\times
f_2] = \Delta'^* [X_1\times X_2, f_1\times f_2]\in \mathcal U^{d_1
+ d_2} (Y)
\]
where $\Delta'$ is a smooth $\varepsilon$-approximation of
$\Delta$ which is transverse to $f_1\times f_2$.
\end{thm}
The cup product is well-defined and associative.

Then, $\mathcal U^* (\,)$ is a multiplicative contravariant
functor for smooth functions on the separable Hilbert manifolds.
The question of whether it agrees with other cobordism functors
such as representable cobordism seems not so easily answered and
there is also no obvious dual bordism functor.

In this section, we show how to define Euler classes in complex
cobordism for finite dimensional complex vector bundles over
separable Hilbert manifolds. In order to do this, we use Sard
classes. We know from \cite{janich} that global sections of a
vector bundle on a smooth separable Hilbert manifold can be
constructed using partitions of unity, then all sections are Sard.
Given a smooth vector bundle $\pi: E\rightarrow B$ over a
separable Hilbert manifold $B$, we know from Theorem \ref{open
embedding}, that $B$ can be embedded as a open subset of a
separable Hilbert space $H$. By Theorem \ref{hist}, we have
\begin{cor}
Let $\pi: E\rightarrow B$ be a finite dimensional complex vector
bundle over a separable Hilbert manifold $B$ and let $i:
B\rightarrow E$ be the zero-section. Then there is an
approximation $\tilde{i}$ of $i$ with $\tilde{i}$ transverse to
$i$.
\end{cor}
Then, we define the Euler class of a finite dimensional complex
vector bundle on a separable Hilbert manifold. Note that Theorem
\ref{son} implies that this Euler class is a well-defined
invariant of the bundle $\pi$.
\begin{defn}
Let $\pi: \xi \rightarrow B$ be a finite dimensional complex
vector bundle of dimension $d$ on a separable Hilbert manifold $B$
with zero-section $i: B\rightarrow \xi$. The \emph{$\mathcal
U$-theory Euler class} of $\xi$ is the element
$$
\chi(\pi) =i^* i_* (1)\in \mathcal U^{2d}(B).
$$
\end{defn}
We have the following projection formula for the Gysin map.
\begin{thm}\label{projection formula}
Let $f: X\rightarrow Y$ be an admissible complex orientable map
and let $\pi: \xi\rightarrow Y$ be a finite dimensional smooth
complex vector bundle of dimension d. Then
\[
\chi (\xi) \cup [X, f] = f_* \chi (f^* \xi).
\]
\end{thm}
\begin{proof}
Let $s$ be a smooth section of $\pi$ transverse to the zero
section $i: Y\rightarrow \xi$. Then $Y'= \{y \in Y: s(y) = i(y)\}$
is a submanifold of complex codimension $d$ and $\chi(\xi) =
[Y',j]$, where $j: Y'\rightarrow Y$ is the inclusion. Setting
$$
X' = f^{-1} Y' = \{x\in X : s(f(x))= i(f(x))\},
$$
which is also a submanifold of $X$ of complex codimension $d$, we
have
\begin{align}
\chi (\xi) \cup [X,f] &= [Y',j] \cup [X,f]\notag\\
& = [X',f_{|X'}].\notag
\end{align}
Now we determine  $f_* \chi (f^* E)$. By transversality theorem,
$f$ can be deformed to a smooth map $f'$ such that the composite
section $s\circ f':X\rightarrow f'^*$ is transverse to the zero
section and they agree on $X'$, hence by definition we have $\chi
(f^* \xi) = [X',j]$ where $j:X'\rightarrow X$ is the inclusion.
Hence, $f_*\chi (f^* \xi) = [X',f_{|X'}]$ by definition of the
Gysin map $f_*$.
\end{proof}
Now we need a useful lemma from \cite{Quinn}.
\begin{lem}\label{tub}
A smooth split submanifold of a smooth separable Hilbert manifold
has a smooth tubular neighborhood.
\end{lem}

Let $\pi: \xi \rightarrow X$ be a finite dimensional complex
vector bundle of dimension $d$ on a separable Hilbert manifold $X$
with zero-section $i: X\rightarrow \xi$. The map $i$ is proper so
that we have the Gysin map
$$
i_*: \mathcal U^{j}(X)\rightarrow \mathcal U^{j+2d}(\xi,\xi-U)
$$
where $U$ is a smooth neighborhood of the zero section.

The map $\pi$ is not proper. However if $U$ is contained in a tube
$U^r$ of finite radius $r$, then $\pi_{|\bar{U}}$ is proper and we
can define
$$
\pi_*:\mathcal U^{j+2d}(\xi,\xi-U)\rightarrow \mathcal U^{j}(X).
$$
Since $\pi i= \Id$ we have $\pi_* i_* = \Id$. The composite map
$i\pi $ is homotopic to $\Id_{\xi}$. If $U = U^{\circ}$ is itself
a tube, the homotopy moves on $U$ and we have \emph{Thom
isomorphism}
$$
\mathcal U^{j+2d}(\xi,\xi-U)\cong \mathcal U^{j}(X).
$$

\section{The relationship between $\mathcal U$-theory and
$MU$-theory.} In this section we consider the relationship between
$\mathcal U$-theory and $MU$-theory. Later we discuss the
particular cases of Grassmannians and $LG/T$.
\par
First we discuss the general relationship between $\U^*(\ )$ and
$MU^*(\ )$. Let $X$ be a separable Hilbert manifold. For each
proper smooth map $f: M\lra X$ where $M$ is a finite dimensional
manifold, there is a pullback homomorphism $$f^*:
\U^*(X)\lra\U^*(M)=MU^*(M).$$ If we consider all such maps into
$X$, then there is a unique homomorphism
$$
\rho :\U^*(X)\lra\ds\lim_{\overleftarrow{M\downarrow X}}MU^*(M),
$$
where the limit is taken over all proper smooth maps $M\rightarrow
X$ from finite dimensional manifolds, which form a directed system
along commuting diagrams of the form
\[
\begin{CD}
M_1  @>f>> M_2  \\
@VVV         @VVV \\
X @>=>> X
\end{CD}
\]
and hence give rise to an inverse system along induced maps $f^* :
MU^* (M_2)\rightarrow MU^* (M_1)$ in cobordism.

Let $X$ be a separable Hilbert manifold. Each of the following
conjectures appears reasonable and is consistent with examples we
will discuss later. We might also hope that surjectivity could be
replaced by isomorphism, but we do not have any examples
supporting this.
\begin{conj}
\label{conj:InvLimSurj} $\rho$ is always a surjection.
\end{conj}
\begin{conj}
If $\U^{\ev}(X)=0$ or $\U^{\odd}(X)=0$, $\rho$ is a surjection.
\end{conj}
\begin{conj}
If $MU^{\ev}(X)=0$ or $MU^{\odd}(X)=0$, $\rho$ is a surjection.
\end{conj}

Now we discuss some important special cases. Let $H$ be a
separable complex Hilbert space, with $H^n$ ($n\geq1$) an
increasing sequence of finite dimensional subspaces with $\dim H^n
= n$ with $H^\infty=\ds\bigcup_{n\geq 1} H^n$ dense in $H$. We use
a theorem of Kuiper \cite{Kuiper}.
\begin{thm}\label{Kuiper}
The unitary group $U(H)$ of a separable Hilbert space $H$ is
contractible.
\end{thm}

Let $\Gr_n(H)$ be the space of all $n$-dimensional subspaces of
$H$, which is a separable Hilbert manifold. Then
$$
\Gr_n(H^\infty) = \ds\bigcup_{k\geq n} \Gr_n (H^k)
$$
is a dense subspace of $\Gr_n (H)$ which we will take it to be a
model for the classifying space $BU(n)$.
\begin{thm}
\label{ProjHilbSpace} The natural embedding $\Gr_n
(H^\infty)\rightarrow \Gr_{n}(H)$ is a homotopy equivalence, and
the natural $n$-plane bundle $\xi_n\rightarrow \Gr_n(H)$ is
universal.
\end{thm}
\begin{proof}
By a theorem of Pressley and Segal \cite{segal-pressley}, the
unitary group $U(H)$ acts on $\Gr(H)$ transitively and hence
$U(H)$ acts on $\Gr_n (H)$ transitively. Let $H^n$ be an
$n$-dimensional subspace of infinite dimensional separable Hilbert
space $H$ and let $H'$ be the orthogonal complement of $H^n$ in
$H$. The stabilizer group of $H^n$ is $U(H^n) \times U(H')$ which
acts freely on the contractible space $U(H)$. Hence
\begin{align}
\Gr_{n}(H)& = U(H)/(U(H^n) \times U(H'))\notag\\
&= B(U(H^n) \times U(H'))\notag\\
&= BU(H^n) \times BU(H').\notag
\end{align}
By Kuiper's Theorem \ref{Kuiper}, $U(H')$ is contractible, hence
so is $BU(H')$. Hence
$$
\Gr_{n} (H) \simeq BU(H^n) = BU(n).
$$
On the other hand,
\[
\Gr_n (H^\infty) = \ds\bigcup_{k\geq n} U(H^k)/(U(H^n) \times
U(H''))\subset \Gr_n(H),
\]
where $H''$ is the orthogonal complement of $H^n$ in $H^k$.

By the construction, the natural $n$-plane bundle
$\xi_n\rightarrow \Gr_n(H)$ is universal. Also, the natural bundle
$\xi_n^\infty\rightarrow \Gr_n(H^\infty)$ is classified by the
inclusion $\Gr_n(H^\infty)\rightarrow \Gr_n(H)$ and since the
latter is universal, this inclusion is a homotopy equivalence.
\end{proof}

In particular, the inclusion of the projective space
$$
P(H^\infty)=\ds\bigcup_{n\geq 1} P(H^n) \subset P(H)
$$
is a homotopy equivalence.
\begin{thm}
The natural homomorphism
$$
\rho: \U^*(P(H))\lra\ds\lim_{\overleftarrow{n}}
MU^*(P(H^n))=MU^*(P (H^\infty))
$$
is surjective.
\end{thm}
\begin{proof}
We will show by induction that
$$
\U^*(P(H))\xrightarrow{i_n^*}MU^*(P(H^{n+1}))
$$
is surjective for each $n$. It will suffice to show that $x^i\in
\im i_n^*$ for $i=0,\ldots,n$. For $n=0$, this is trivial.

Now we verify it for $n=1$. By Theorem \ref{ProjHilbSpace}, since
the natural line bundle $\lambda\rightarrow P(H)\simeq
P(H^\infty)$ is universal, the following diagram commutes for each
$n\geq 1$
\[
\begin{CD}
\eta_n =i_n^*(\lambda) @>i_n^*>> \lambda\\
@VVi_n^*(\lambda)V  @VVV\\
\mathbb C P^n = P(H^{n+1}) @>i_n>> P(H),
\end{CD}
\]
where $i_n: \mathbb C P^{n} =P(H^{n+1}) \rightarrow P(H)$ is the
inclusion map. By the compatibility of induced bundles, for $n\geq
1$ and the generator $x = \chi(\eta_n)\in MU^*(P(H^{n+1}))$, there
exists an Euler class $\tilde{x} = \chi(\lambda) \in \mathcal
U^2(P(H))$ satisfying $i_n^*(\tilde{x}) = x$, where $i_n:
P(H^{n+1}) \rightarrow P(H)$ is the inclusion map.

Assume that $i_n^*$ is surjective. Then there are elements
$$
y_i\in \U^{2i}(P(H)),\, i=0,\ldots,n,
$$
such that
$$
i_n^*y_i = x^i \in MU^{2i} (P(H^{n+1})).
$$
Also,
$$
i_{n+1}^*  y_i = x^i + z_i  x^{n+1} \in
MU^{2i}(P(H^{n+2}))
$$
where $z_i \in MU_{2(n+1-i)}$.

In particular, let $y_{n} = [W, f] \in \U^{2n}(P(H))$. Then the
following diagram commutes
\[
\begin{CD}
f^*(\lambda) @>f^*>> \lambda\\
@VVV  @VVV\\
W @>f>> P(H)
\end{CD}
\]
and there is an Euler class $\chi(f^*(\lambda)) = [W', g] \in
\U^2(W)$. Now by Theorem \ref{projection formula},
$$
y_{n+1} = f_* \chi(f^*(\lambda)) \in \U^{2n+2}(P(H))
$$
satisfies
$$
i_{n+1}^* y_{n+1} = x^{n} \chi (\eta_n) =x^{n+1}.
$$

Hence, $\im i_{n+1}^*$ contains the $MU^*$-submodule generated by
$x^i$ $(i=0,\ldots,n)$ and so $i_{n+1}^*$ is surjective. This
completes the induction.

This shows that the induced homomorphism
$$
\rho:\U^*(P(H))\lra\ds\lim_{\overleftarrow{n}} MU^*(P(H^n))=MU^*(P
(H^\infty))
$$
is surjective.
\end{proof}

Note that it is also possible to prove this result by using the
projective spaces $P({H^n}^{\bot}) \subset P(H)$ to realize
cobordism classes restricting to the classes $x^n$ in $MU^*(P
(H^\infty))$.

Next we discuss some geometry of Grassmannians from Pressley $\&$
Segal \cite{segal-pressley}, whose ideas and notation we assume.
We take for our separable Hilbert space $H = L^2 (\S^1 ; \C)$ and
let $H_+$ to be the closure of the subspace of $H$ containing the
functions $z^n : z\rightarrow z^n (n\geq 0)$.  Then
$$
\Gr_0 (H) = \ds\lim_{\overrightarrow{k\geq 1}} \Gr (H_{-k,k}),
$$
where $\Gr (H_{-k,k})$ is the Grassmannian of the finite
dimensional vector space
$$
H_{-k,k} = z^{-k}H_+ / z^k H_+.
$$
$\Gr_0 (H)$ is dense in $\Gr (H)$ and is also known to be
homotopic to the classifying space of $K$-theory, $BU\times \Z$.
\begin{thm}
For $n\geq 1$, the natural homomorphism
$$
\rho: \mathcal
U^*(\Gr_{n}(H))\rightarrow MU^*(Gr_n(H))
$$
is surjective.
\end{thm}
\begin{proof}
For $k\geq n$, the inclusion $i: \Gr_n(H_{-k,k})\rightarrow
\Gr_n(H)$ induces a contravariant map
$$\U^*(Gr_n(H)) \rightarrow \U^*(\Gr_n(H_{-k,k}) =MU^*(\Gr_n(H_{-k,k})).$$
For $k\geq n$, since $C_S \subset \Gr_n(H_{-k,k})$ is transverse
to $\Sigma_S$, there exists a stratum $\Sigma_{S'}$ such that
\[
\sigma_{S',k}= [\Gr_n (H_{-k,k}) \cap \Sigma_{S'} \rightarrow
\Gr_n(H_{-k,k})] \in MU^*(\Gr_n(H_{-k,k}))
\]
are the classical Schubert cells. By an argument using the
Atiyah-Hirzebruch spectral sequence and results on Schubert cells
in cohomology \cite{milnor-stasheff}, the cobordism classes
$\sigma_{S', k}$ provide generators for the $MU^*$-module
$MU^*(\Gr_n(H_{-k,k}))$. Then $i^*$ is surjective. For each $k$,
$$MU^{\mathrm{odd}}(\Gr_n(H_{-k,k}))= 0,$$ hence
\begin{align}
\U^*(\Gr_n(H))\rightarrow& \ds\lim_{\overleftarrow{k}}MU^*(\Gr_n(H_{-k,k}))\notag\\
&= MU^*(\Gr_n(H^\infty))\notag\\
&\cong MU^*(\Gr_n(H))\notag
\end{align}
is surjective.
\end{proof}

\begin{thm}
For a compact connected semi-simple Lie group $G$,
$$
\rho: \U^*(LG/T)\rightarrow MU^*(LG / T)
$$
is surjective.
\end{thm}
\begin{proof}
As $LG/T$ has no odd dimensional cells,the Atiyah- Hirzebruch
spectral sequence for  $MU^*(LG/T)$ collapses. Hence it suffices
to show that the  composition
\[
\U^*(LG/T) \rightarrow MU^*(LG/T) \rightarrow H^*(LG/T, \mathbb Z)
\]
is surjective. Since $H^*(LG/T,\mathbb Z)$ is generated by the
Schubert classes $\varepsilon^w (w \in W)$ dual to the Schubert
cells $C_w$, and $\Sigma_w$ is dual to $C_w$, the image of the
stratum $\Sigma_w$ under the composition map gives
$\varepsilon^w$, establishing the desired surjectivity.
\end{proof}

Similarly, we have
\begin{thm}
For a compact connected semi-simple Lie group $G$,
$$
\rho : \U^*(\Omega G) \rightarrow MU^*(\Omega G)
$$
is surjective.
\end{thm}

\section{Cobordism classes related to Pressley-Segal
stratifications and some calculations.} In this section, we show
that the stratifications introduced by Pressley $\&$ Segal
\cite{segal-pressley} give rise some further interesting cobordism
classes in $\mathcal U^*(LG /T)$.
\par
The Grassmannian $\Gr H$ of \cite{segal-pressley} is a separable
Hilbert manifold, and the stratum $\Sigma_S \subset \Gr H$ is a
locally closed contractible complex submanifold of codimension
$\ell (S)$ and the inclusion map $\Sigma_S \rightarrow \Gr(H)$ is
a proper Fredholm map with index $- \ell (S)$.
\par
Therefore, we have
\begin{thm}
The stratum $\Sigma_S\rightarrow \Gr (H)$ represents a class in
$\mathcal U^{2 \ell (S)} (\Gr (H))$.
\end{thm}
These strata $\Sigma_S$ are dual to the Schubert cells $C_S$ in
the following sense;
\begin{enumerate}
\item the dimension of $C_S$ is the codimension of
$\Sigma_S$ and
\item $C_S$ meets $\Sigma_S$ transversely in a single point, and
meets no other stratum the same codimension.
\end{enumerate}
\par
The loop group $LG$ acts via the adjoint representation on the
Hilbert space
$$
H_{\mathfrak g} = L^2 (\S^1 ; {\mathfrak g}_{\C}),
$$
where ${\mathfrak g}_{\C})$ is the complexified Lie algebra of
$G$. If $\dim G = n$, we can identify $H_{\mathfrak g}$ with $H^n$
and since the adjoint representation is unitary for a suitable
Hermitian inner product, this identifies $LG$ with a subgroup of
$LU(n)$. Then \cite{segal-pressley} shows how to identify the
based loop group $\Omega G$ with a submanifold of $\Omega U(n)$,
which can be itself identified with a submanifold of $\Gr
(H_{\mathfrak g})$.

Then $\Omega G$ inherits a stratification with strata
$\Sigma_{\lambda}$ indexed by homomorphisms $\lambda :
\T\rightarrow T$. Each stratum $\Sigma_{\lambda} \subset  \Omega
G$ is a locally closed contractible complex submanifold of
codimension $d_{\lambda}$, and the inclusion map
$\Sigma_{\lambda}\rightarrow \Omega G$ is an admissible Fredholm
map. Then
\begin{thm}
For each $\lambda$, the inclusion $\Sigma_{\lambda}\rightarrow
\Omega G$ represents a class in $\mathcal U^{2 d_{\lambda}}
(\Omega G)$.
\end{thm}

If we restrict to the inverse limit $MU^*(\Omega G)$, these
stratum $\Sigma_{\lambda}$ provide the basis elements of
$MU^{2d_{\lambda}}(\Omega G)$. Since $H^*(\Omega G)$ is even
graded and $MU^*$ is also even graded, the Atiyah-Hirzebruch
spectral sequence
$$
H^*(\Omega G; MU^*)\Rightarrow MU^*(\Omega G)
$$
collapses, hence we have an isomorphism
$$
\mathcal U^* (\Omega G) /{\ker \rho} \cong MU^*(\Omega G) \cong
H^*(\Omega G)\otimes MU^*.
$$
For $G= SU(2)$, we have
$$
\mathcal U^* (\Omega SU(2)) /{\ker \rho} \cong MU^*(\Omega SU(2))
\cong \Gamma_{\Z}(\gamma)\otimes MU^*,
$$
where $\Gamma_{\Z}(\gamma)$ is a divided power algebra with the
$\Z$-module basis $\gamma^{[n]}$ in each degree $2n$ for $n\geq
1$.

\par
Such stratifications also exist for the homogeneous space $LG /T$.
\begin{thm}For $w \in
\widetilde{W}$, the inclusion $\Sigma_{w}\rightarrow LG/T$
represents a class in $\mathcal U^{2 \ell (w)} (L G / T)$.
\end{thm}

Similarly, if we restrict to the inverse limit $MU^*(LG/T)$, we
have an isomorphism
$$
\mathcal U^* (LG/T) /{\ker \rho} \cong MU^*(LG/T) \cong H^*(
LG/T)\otimes MU^*.
$$
For $G= SU(2)$, we have
$$
\mathcal U^* (LSU(2)/T) /{\ker \rho} \cong MU^*(LSU(2)/T) \cong
\Gamma_{\Z}(x_0, x_1)/I_{\Z}\otimes MU^*,
$$
where the ideal $I_{\Z}$ is given by
$$
I_{\Z} = \left(x_0^{[n]} x_1^{[m]} - \binom{n + m - 1}{m} x_0^{[n
+ m]} -  \binom{n + m - 1}{n} x_1^{[n + m]}: \, m,n \geq 1\right)
$$
and which has the $\Z$-module basis $\{x_0^{[n]}, x_1^{[n]}\}$ in
each degree $2 n$ for $n \geq 1$.

\end{document}